# Accurate Semidefinite Programming Models for Optimal Power Flow in Distribution Systems

Zeyu Wang, Daniel S. Kirschen, *Fellow, IEEE*, Baosen Zhang

***Abstract*****-- In this paper, we develop semidefinite programming (SDP) models aimed at solving optimal power flow (OPF) problems in distribution systems. We first propose a symmetrical SDP model which modifies the existing BFM-SDP model. Comprehensive case studies are conducted to show that our SDP approach is more numerically stable and more accurate than existing approaches. Based on the symmetrical SDP model, we also develop a technique to solve the problem of keeping the nodal voltages within their bounds. By comparing our results to benchmark power flow solutions generated using OpenDSS, we provide a rigorous assessment of our claims.**

## I. Introduction

As the penetration of distributed energy resources increases, distribution networks must be actively managed to ensure their reliability and optimize their operation. Techniques that have been used to solve the optimal power flow (OPF) problem in transmission networks are not directly applicable to distribution networks because of three fundamental differences between these networks: 1) most distribution networks are operated radially; 2) the three phases are usually unbalanced in distribution networks 3) loads at the distribution level cannot be modeled as being independent of the voltage. To the authors' knowledge, no full-fledged OPF solver can deliver accurate results on unbalanced three-phase distribution systems. This paper proposes a more accurate and numerically stable semi-definite programming (SDP) OPF solution technique that specifically addresses these differences. In particular, this technique:

- Uses symmetrical components to reduce the coupling between the phases in the backbones of the distribution network.
- Models the tap-changing voltage regulators to account for their effect on the voltage profile.
- Iteratively updates the loads to reflect their dependence on the voltage
- Accurately handles binding voltage constraints.

The remainder of this paper is organized as follows. Section II reviews the methods that have been proposed to solve the OPF problem in distribution networks. Section III describes the proposed symmetrical SDP formulation and contrasts it with existing methods. Section IV demonstrates that the symmetrical SDP method is more accurate and numerically stable than other techniques using comprehensive case studies on several IEEE test feeders. Section V describes a solution to the voltage regulation problem. Section VI validates this approach using case studies. Section VII draws conclusions and suggests directions for further research.

## II. Literature Review

Several methods have been proposed to overcome the non-linear and non-convex nature of the OPF problem in distribution networks. Baran and Wu [1] develop a linearized OPF model similar to the DC power flow in transmission systems. Jabr [2] proposed a conic model to solve distribution system power flow. Based on this approach, Farivar and Low [3,4] developed the branch flow model, a second order cone programming (SOCP) OPF model. Gan et. al [5] introduced an SOCP model for multiphase distribution networks. However, this model decomposes the three-phase network into three independent single-phase networks and thus ignores the coupling between phases, possibly leading to inaccurate results.

Three-phase OPF models have been developed to better represent distribution systems. Bruno et. al [6] proposed an iterative Newton method to adjust the decision variables and utilizes a power flow solver to update the state variables. Dall'Anese et. al [7] applied semidefinite programming (SDP) on a bus-injection model (BIM-SDP). Gan and Low [8] enhanced this model by applying chordal relaxation which reduces the number of variables and speeds up the solution process. Later on, Gan and Low [8] proposed a branch-flow model (BFM-SDP) that is equivalent to BIM-SDP but has better numerical stability. They also developed a linearized three-phase OPF model (LPF) which extends Baran and Wu's model [1] to three phase networks. However, this model ignores losses, which can be as high as 10% in distribution networks. Zamzam et al. [9] developed a QCQP model that replaces the non-convex part of the OPF constraints by convex approximations. However, as these authors observe, without warm-start, the method could take 1000 iterations to converge. Table I summarizes the characteristics of the OPF methods described above.

Table I. Summary of Distribution Network OPF Methods

| | **Model** | **Phase** | **Losses** | **Comments** |
|---|---|---|---|---|
| DISTFLOW [1] | Linear | Single | No | Single phase, linear |
| SOCP [2,3,4,5] | SOCP | Single /Three* | Yes | *Break three phase down to 3 single phase models, essentially single phase |
| BIM-SDP [7,8] | SDP | Three | Yes | Numerical issues |
| BFM-SDP [8] | SDP | Three | Yes | More stable than BIM-SDP |
| 3ϕ-linear LPF-OPF [8] | Linear | Three | No | Ignores losses |
| Iterative [6] | NR | Three | Yes | Requires PF solver |
| Successive convex approximation[9] | QCQP | Three | Yes | Very slow without warm start |

Model: category of the convex optimization technique; Phase: single-phase or three-phase model; Losses: whether this method considers losses.

It is also useful to relate these techniques to methods that have been proposed recently to solve the OPF problem in transmission systems. Sandro et al. [10] applied an approach similar to the ideas expounded in [9] on transmission systems. Molzahn and Hiskens [11] developed a moment-based relaxation that is tighter than conventional SDP programs, but at the computational cost of larger semidefinite programs. Coffrin et al. [12] proposed a QC relaxation to solve OPF problems. This reference also presents a comprehensive evaluation of the performance of OPF techniques in transmission systems.

To assess the quality of the solution produced by these various OPF methods, their solutions (i.e. the dispatch of distributed energy resources) can be input to the well-known distribution system power flow solver OpenDSS [13]. The results from OpenDSS power flow calculations thus provide a benchmark to assess the accuracy of OPF solutions. Emiroglu et. al [14] suggested that linear models produce solutions that are not accurate, especially with regard to voltage magnitudes. While the BFM-SDP model is more stable, section IV shows its results differ from the OpenDSS solution in most test cases. This paper proposes a method that yield very accurate results.

## III. Symmetrical SDP model formulation

This section first reviews briefly the BFM-SDP model, then describes the three essential components of the proposed symmetrical SDP technique: the use of symmetrical components, the inclusion of voltage regulators and an accurate modeling of loads.

### A. BFM-SDP Model

The BFM-SDP model represent multiphase voltages and currents as vectors. For each line segment from node *i* to node *j*, we define the nodal voltage vectors $V_i$ and $V_j$ , and the current $I_{ij}$ . For example, if node *i* has three phases, then $V_i=[V_i^a,V_i^b,V_i^c]^T$ . If node *j* only has phases a and c, then $V_j=[V_j^a,V_j^c]^T$ and *ij* is a two-phase line: $I_{ij}=[I_{ij}^a,I_{ij}^c]^T$ .

We define the second-order decision variables using matrices: $v_i=V_iV_i^H$ , $\ell_{ij}=I_{ij}I_{ij}^H$ and $S_{ij}=V_i^{\Phi_{ij}}I_{ij}^H$ . In the above example, $v_i$ is a 3x3 matrix, while $\ell_{ij}$ and $S_{ij}$ are 2x2 matrices. $\Phi_{ij}$ is the set of phases of line *ij*. The superscript *H* indicates the Hermitian transpose. The branch flow SDP model can be formulated as:

$$\min \sum_{i\in N} C_i(s_i)$$

$$over: s_i\in\mathbb{C}^{|\Phi_i|}, v_i\in\mathbb{H}^{|\Phi_i|\times|\Phi_i|} \text{ for } i\in N \quad (1)$$

$$S_{ij}\in\mathbb{C}^{|\Phi_{ij}|\times|\Phi_{ij}|}, \ell_{ij}\in\mathbb{H}^{|\Phi_{ij}|\times|\Phi_{ij}|} \text{ for } i\to j$$

where $\mathbb{C}$ is the set of complex numbers and $\mathbb{H}$ is the set of Hermitian matrices. $s_i$ is the nodal injection at node *i*, i.e. either loads or the net injection of a distributed energy resource. This optimization is subject to the following constraints:

$$\sum_{i:i\to j} diag(S_{ij}-z_{ij}l_{ij})+s_j+y_jv_j=\sum_{k:j\to k} diag(S_{jk})^{\Phi_j} \ j\in N \quad (2)$$

Where Equation (2) is the power flow balance constraint. The first term represents the upstream power flow and line losses along line *ij*, the last term represents the downstream power flows from node *j*. The operation $diag(A)$ returns the diagonal vector of matrix $A$ . $z_{ij}$ is the line impedance and $y_j$ the nodal shunt capacitance. Equation (3) describes Kirchoff's voltage law along line *ij*:

$$v_j=v_i^{\Phi_{ij}}-(S_{ij}z_{ij}^H+S_{ij}^Hz_{ij})+z_{ij}\ell_{ij}z_{ij}^H \quad i\to j \quad (3)$$

Equation (4) enforces the lower and upper bounds on the nodal voltage while Equation (5) sets the voltage at the source node.

$$\underline{v}_i\le diag(v_i)\le\overline{v}^i \ \ i\in N \quad (4)$$

$$v_0=V_0^{ref}(V_0^{ref})^H \quad (5)$$

Equation (6) is the positive semidefinite constraint. Equation (7) enforces that this positive semidefinite matrix should be of rank one.

$$\begin{bmatrix} v_i^{\Phi_{ij}} & S_{ij} \\ S_{ij}^H & \ell_{ij} \end{bmatrix}\succeq 0 \ i\to j \quad (6)$$

$$\text{rank}\begin{bmatrix} v_i^{\Phi_{ij}} & S_{ij} \\ S_{ij}^H & \ell_{ij} \end{bmatrix}=1 \ i\to j \quad (7)$$

The above formulation is equivalent to the original, non-convex OPF formulation described in [8]. The BFM-SDP approach removes the non-convex constraint (7) and solves the SDP model with the convex objective function $C_i(s_i)$ and the convex constraints (2-6). If the rank-1 constraint (7) holds for a given SDP solutions, the SDP relaxation is tight, and the solution is equivalent to the original OPF solution. On the other hand, if (7) doesn't hold the SDP solution does not have a direct physical meaning.

### B. Symmetrical SDP Model

The symmetrical components transformation reduces the phase coupling in three-phase systems and has been used for decades to study unbalanced systems [15]. It has been applied to the solution of the unbalanced three-phase power flow problem [16]. In this section, we show how its application to the three-phase backbone of distribution networks enhances the numerical stability of an SDP-based OPF method.

Voltages in phase components are transformed into symmetrical components as follows:

$$V^{abc}=AV^{012} \quad (8)$$

Where:

$$A = \frac{1}{\sqrt{3}}\begin{bmatrix} 1 & 1 & 1 \\ 1 & a^2 & a \\ 1 & a & a^2 \end{bmatrix}, \ A^H = A^{-1} \tag{9}$$

and $a = 1\angle 120°$. Hence, the three-phase variables in the BFM-SDP method and the impedance parameters are related to the equivalent variables in symmetrical components:

$$\begin{aligned} v^{abc} &= V^{abc} \times V^{abc,H} \\ &= AV^{012} \times (AV^{012})^H \\ &= Av^{012}A^H \end{aligned} \tag{10}$$

$$\begin{aligned} l^{abc} &= Al^{012}A^H \\ S^{abc} &= AS^{012}A^H \\ z^{abc} &= Az^{012}A^{-1} = Az^{012}A^H \\ y^{abc} &= Ay^{012}A^{-1} = Ay^{012}A^H \end{aligned} \tag{11}$$

Constraints (2) – (6) are transformed as follows (12) – (16):

$$\begin{aligned} &\sum_{i:i\to j} diag(A(S_{ij}^{012} - z_{ij}^{012}l_{ij}^{012})A^H) + s_j + y_j^{012}v_j^{012} \\ &= \sum_{k:j\to k} diag(AS_{jk}^{012}A^H) \end{aligned} \tag{12}$$

$$v_j^{012} = v_i^{012} - (S_{ij}^{012}z_{ij}^{012,H} + S_{ij}^{012,H}z_{ij}^{012}) + z_{ij}^{012}\ell_{ij}^{012}z_{ij}^{012,H} \tag{13}$$

$$\underline{v}_i \le diag(Av_i^{012}A^H) \le \overline{v}^i \tag{14}$$

$$v_0^{012} = V_0^{012,ref}(V_0^{012,ref})^H \tag{15}$$

$$\begin{bmatrix} v_i^{012} & S_{ij}^{012} \\ S_{ij}^{012,H} & \ell_{ij}^{012} \end{bmatrix} \succeq 0 \tag{16}$$

For three-phase nodes that connect to single- or two-phase laterals, we introduce auxiliary variables that convert symmetrical component variables back to phase variables:

$$v_i^{abc} = Av_i^{012}A^H \tag{17}$$

These auxiliary variables define the boundaries of the symmetrical component backbones. For two-phase and single-phase laterals, equations (2) – (6) still apply. The symmetrical SDP OPF model thus consists of the objective function (1), subject to constraints (2) – (6) on two- and single-phase laterals and constraints (12) – (16) on the three-phase backbones.

### C. Modeling Voltage regulators

American National Standards Institute (ANSI) standard ANSI C84.1 [17] defines the range of acceptable distribution voltages in the United States. Step-voltage regulators, which are essentially tap-changing transformers, are commonly used to maintain system voltages within this range. For three single-phase voltage regulators installed at a three-phase bus, the ratios between the primary and secondary voltages are given by:

$$\begin{aligned} &ratio = [r_a, r_b, r_c]^T \\ &[V_{\sec}^a, V_{\sec}^b, V_{\sec}^c]^T = [r_aV_{pri}^a, r_bV_{pri}^b, r_cV_{pri}^c]^T \end{aligned} \tag{18}$$

where

$$\begin{aligned} r_a &= 1 + 0.00625 \times Tap_a \\ r_b &= 1 + 0.00625 \times Tap_b \\ r_c &= 1 + 0.00625 \times Tap_c \end{aligned} \tag{19}$$

$Tap_a$, $Tap_b$ and $Tap_c$ are integers between $[-16,16]$.

The voltages on the two sides of a regulator are linked as follows:

$$v_{\sec}^{abc} = (v_i^{\Phi_{ij}} - (S_{reg}z_{reg}^H + S_{reg}^Hz_{reg}) + z_{reg}\ell_{reg}z_{reg}^H) \bullet \mathrm{R} \tag{20}$$

where $\bullet$ indicates a dot product and $\mathrm{R} = ratio \times ratio^T$ is a 3x3 matrix and $z_{reg}$ is the impedance of the voltage regulator. In terms of symmetrical components, (20) becomes:

$$\begin{aligned} &Av_{\sec}^{012}A^H = \\ &A[v_{pri}^{012} - (S_{reg}^{012}z_{reg}^{012,H} + S_{reg}^{012,H}z_{reg}^{012}) + z_{reg}^{012}\ell_{reg}^{012}z_{reg}^{012,H}]A^H \bullet R \end{aligned} \tag{21}$$

At voltage regulator nodes, we replace the voltage constraints (2) and (13) by (20) and (21) to integrate the three-phase voltage regulators into the model. A similar approach is used for single- and two-phase voltage regulators.

### D. Modeling Loads

Power flow calculations typically assume that loads can be modeled as requiring a given amount of active and reactive power. In transmission networks this is a good approximation because the voltage at the nodes where these loads are modeled is held constant through the action of tap-changing transformers. On the other hand, this approximation is not valid within a distribution network because the voltages at the nodes where these loads are located cannot be assumed constant. Models of loads should therefore reflect their constant impedance, constant current and constant power characteristics (a.k.a. ZIP models). In addition, wye- and delta-connected loads must also be handled correctly. Figure 1 illustrates the iterative procedure used update the values of loads based on the voltages calculated by the OPF results at the previous iteration.

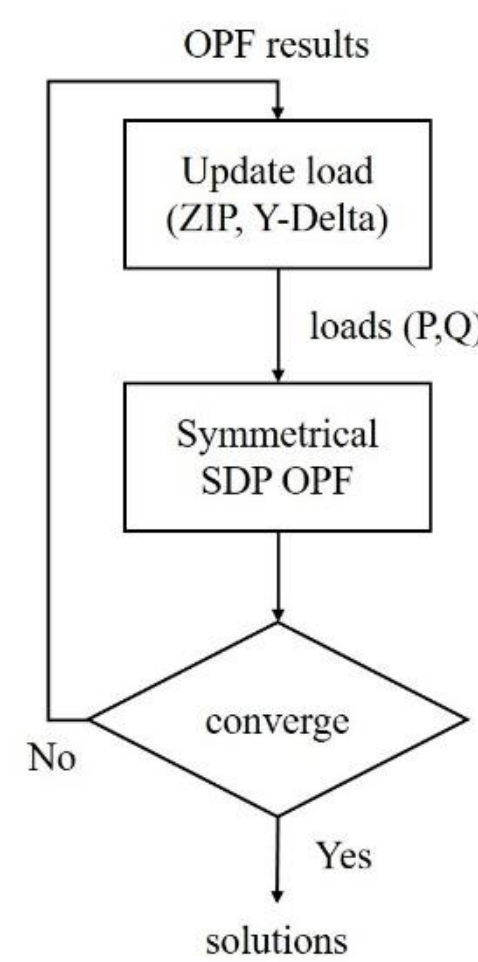

Figure 1. Flow chart of the iterative load update procedure.

## IV. CASE STUDIES

Case studies were performed on the IEEE 13-, 34-, 37-, 123- and 8500-node test feeders [18]. Several modifications were made to these feeder models:

- The delta-connected voltage regulator of the IEEE 37-node test feeder was replaced by a 100ft long cable with line code 721
- Only the primary circuits of the 8500-node test feeder, was modeled. Loads connected to the secondary were directly connect to primary nodes

- The OpenDSS load parameter "Vminpu" was set to 0.85 to prevent the load model reverting to constant impedance
- The tap positions of the voltage regulators in the SDP models were aligned with the IEEE feeder documents and OpenDSS settings.

Calculations were performed using CVX [19] integrated in MATLAB.

### A. *Voltage Profile*

In a first set of studies, no distributed energy resources were incorporated the networks. Solving an OPF is equivalent to solving a power flow. This allows us to compare the results of the BFM-SDP and the symmetrical SDP against the OpenDSS benchmark. The BFM-SDP and the symmetrical SDP OPF were solved with the same objective function, i.e. minimize total network losses. Loss minimization is chosen as the objective function because by doing so, the solution obeys constraint (7) thus has physical meanings [20]. For this test, we used the 34-node test feeder because it includes long overhead lines that causes severe voltage drops and requires the use of two regulators to boost the voltages. To ensure a fair comparison, voltage regulators were modeled identically in the BFM-SDP and symmetrical SDP. Table III summarizes the per unit voltage magnitudes and voltage phase angles for phase A. We don't present the phase B and C voltages for brevity. These results show that the voltages calculated using the symmetrical SDP are closer to those obtained with OpenDSS, particularly at the bold nodes (first column of Table III). This improved accuracy is primarily due to the iterative load update process. Ten or fewer iterations of the symmetrical SDP were required for all test cases.

Table III. Nodal voltages for the IEEE 34 Node Test Feeder (Phase A)

| | OpenDSS | | Symmetrical SDP | | BFM-SDP [8] | |
|---|---|---|---|---|---|---|
| Node | VA (p.u.) | thetaA (°) | VA (p.u.) | thetaA (°) | VA (p.u.) | thetaA (°) |
| SOURCE | 1.0000 | 0.0 | 1.0000 | 0.0 | 1.0000 | 0.0 |
| 800 | 1.0000 | 0.0 | 0.9999 | 0.0 | 0.9999 | 0.0 |
| 802 | 0.9974 | -0.1 | 0.9974 | -0.1 | 0.9971 | 0.0 |
| 806 | 0.9957 | -0.1 | 0.9956 | -0.1 | 0.9951 | -0.1 |
| 808 | 0.9631 | -0.7 | 0.9631 | -0.7 | 0.9590 | -0.7 |
| 812 | 0.9253 | -1.6 | 0.9253 | -1.6 | 0.9171 | -1.5 |
| 814 | 0.8952 | -2.3 | 0.8953 | -2.3 | 0.8840 | -2.2 |
| 814R | 0.9624 | -2.3 | 0.9622 | -2.3 | 0.9503 | -2.2 |
| 850 | 0.9624 | -2.3 | 0.9622 | -2.3 | 0.9503 | -2.2 |
| 816 | 0.9619 | -2.3 | 0.9618 | -2.3 | 0.9499 | -2.2 |
| 818 | 0.9610 | -2.3 | 0.9608 | -2.3 | 0.9489 | -2.2 |
| 824 | 0.9529 | -2.4 | 0.9527 | -2.4 | 0.9398 | -2.3 |
| 820 | 0.9361 | -2.3 | 0.9360 | -2.4 | 0.9233 | -2.3 |
| 828 | 0.9522 | -2.4 | 0.9520 | -2.5 | 0.9389 | -2.3 |
| 822 | 0.9328 | -2.3 | 0.9327 | -2.4 | 0.9200 | -2.3 |
| 830 | 0.9343 | -2.6 | 0.9341 | -2.7 | 0.9190 | -2.5 |
| 854 | 0.9338 | -2.6 | 0.9337 | -2.7 | 0.9185 | -2.5 |
| 852 | 0.9030 | -3.1 | 0.9029 | -3.1 | 0.8845 | -2.9 |
| **852R** | 0.9764 | -3.1 | 0.9761 | -3.2 | 0.9563 | -2.9 |
| **832** | 0.9764 | -3.1 | 0.9761 | -3.2 | 0.9563 | -2.9 |
| **858** | 0.9741 | -3.1 | 0.9738 | -3.2 | 0.9537 | -3.0 |
| **888** | 0.9400 | -4.7 | 0.9396 | -4.8 | 0.9173 | -4.8 |
| **864** | 0.9741 | -3.1 | 0.9738 | -3.2 | 0.9537 | -3.0 |
| **834** | 0.9714 | -3.2 | 0.9711 | -3.3 | 0.9508 | -3.0 |
| **890** | 0.8571 | -5.3 | 0.8566 | -5.3 | 0.9148 | -4.8 |
| **860** | 0.9710 | -3.2 | 0.9707 | -3.3 | 0.9503 | -3.0 |
| **842** | 0.9714 | -3.2 | 0.9711 | -3.3 | 0.9507 | -3.0 |
| **836** | 0.9708 | -3.2 | 0.9705 | -3.3 | 0.9501 | -3.0 |
| **844** | 0.9712 | -3.2 | 0.9708 | -3.3 | 0.9505 | -3.1 |
| **840** | 0.9707 | -3.2 | 0.9704 | -3.3 | 0.9501 | -3.0 |
| **862** | 0.9708 | -3.2 | 0.9705 | -3.3 | 0.9501 | -3.0 |
| **846** | 0.9714 | -3.3 | 0.9711 | -3.4 | 0.9507 | -3.1 |
| **848** | 0.9714 | -3.3 | 0.9711 | -3.4 | 0.9507 | -3.1 |

### B. *Numerical Stability*

To demonstrate the superior numerical stability of the symmetrical SDP over the BFM-SDP, distributed generation (DG) was inserted in the feeder models. With the SeDuMi [21] open-source SDP solver, for test most cases, the solution status of BFM-SDP model was "*Inaccurate/Solved*". The solver log shows that it runs into numerical problems after several iterations and terminates with a sub-optimal result. No such numerical issues were encountered with the symmetrical SDP method, and all test cases terminated with the solution status "*Solved*". With the commercial SDP solver Mosek [22], the BFM-SDP approach failed to solve on any test feeder. On the other hand, the symmetrical SDP method succeeded in solving all the test feeders with Mosek. Table II gives the solution time of both solvers on various test feeders. It is interesting that Mosek solves the 123- and 8500-node feeders in less than half of the time needed by SeDuMi.

Table II. OPF Solution Time on Various Test Feeders.

| Solver: | SeDuMi | | Mosek | |
|---|---|---|---|---|
| Test Feeder | BFM-SDP | Sym. SDP | BFM-SDP | Sym. SDP |
| 13node | <1s | <1s | X | <1s |
| 34node | <1s | <1s | X | <1s |
| 37node | 1.56s | 1.06s | X | 0.77s |
| 123node | 4.428s | 4.202s | X | 1.33s |
| 8500node | 31.36s | 43.34s | X | 13.67s |

Note: the solution time accounts only for the time the solver consumes to reach the solutions. The time consumes by CVX to pre-compiles the model is not included. To ensure a fair comparison, the table lists the solution time of only one iteration for both methods. The iterative load update procedure (Figure 1) requires several iterations.

### C. *Accuracy of Optimal Power Flow Solutions*

To gauge the accuracy of the OPF solutions, we compare the three-phase complex power flows at the heads of the feeder, i.e. at the sending end of the lines or the primary side of the voltage regulator located immediately downstream from the primary substation. The results obtained with the BSM-SDP OPF and the Symmetrical SDP OPF are compared with the OpenDSS benchmarks based on the absolute percentage errors on the active and reactive power flows:

$$Error_P^{\phi} = \left| \frac{P_{SDP}^{\phi} - P_{OpenDSS}^{\phi}}{P_{OpenDSS}^{\phi}} \right| \times 100\% \tag{22}$$

$$Error_Q^{\phi} = \left| \frac{Q_{SDP}^{\phi} - Q_{OpenDSS}^{\phi}}{Q_{OpenDSS}^{\phi}} \right| \times 100\% \tag{23}$$

Nine cases studies were performed on different test feeders, with and without voltage regulators, and with substation voltage magnitude set at different levels. Table IV summarizes the settings of the 9 case studies:

Table IV. Settings of 9 Cases on Accuracy of Feeder Head Power Transfer

| Case No. | Test Feeder | Substation Voltage Magnitude p.u. | Consider Regulator | Component Segment of Interest |
|---|---|---|---|---|
| 1 | 13-node | 1.05 | No | Line 650-632 |
| 2 | 13-node | 1.00 | Yes | Line 650-632 |
| 3 | 34-node | 1.05 | Yes | Line 800-802 |
| 4 | 34-node | 1.00 | Yes | Line 800-802 |
| 5 | 37-node | 1.05 | No | Regulator 799R |
| 6 | 37-node | 1.00 | No | Regulator 799R |
| 7 | 123-node | 1.05 | No | Regulator 150R |
| 8 | 123-node | 1.00 | No | Regulator 150R |
| 9 | 123-node | 1.00 | Yes | Regulator 150R |

Tables V and VI show that the symmetrical SDP method achieves substantially less error compared with the BFM-SDP method, especially on the test cases numbers that are bold. In particular, symmetrical SDP method achieves an error on the active power of less than 0.2% in most cases.

Table V. Absolute Percentage Errors of Active Powers at Feeder Head

| Case | Branch Flow SDP | | | Symmetrical SDP | | |
|---|---|---|---|---|---|---|
| No. | PhaseA | PhaseB | PhaseC | PhaseA | PhaseB | PhaseC |
| 1 | 0.457% | 0.930% | 2.570% | 0.008% | 0.102% | 0.007% |
| 2 | 0.351% | 1.650% | 2.182% | 0.016% | 0.000% | 0.007% |
| 3 | 1.217% | 2.682% | 2.553% | 0.029% | 0.059% | 0.036% |
| **4** | 6.979% | 3.652% | 5.453% | 0.017% | 0.110% | 0.026% |
| 5 | 0.008% | 2.519% | 1.095% | 0.019% | 0.044% | 0.046% |
| **6** | 2.653% | 1.862% | 0.823% | 0.128% | 0.079% | 0.182% |
| 7 | 0.130% | 2.197% | 0.933% | 0.021% | 0.018% | 0.008% |
| **8** | 3.997% | 1.250% | 2.771% | 0.014% | 0.019% | 0.000% |
| 9 | 0.396% | 1.519% | 0.916% | 0.055% | 0.154% | 0.286% |

Table VI. Absolute Percentage Errors of Reactive Powers at Feeder Head

| Case | Branch Flow SDP | | | Symmetrical SDP | | |
|---|---|---|---|---|---|---|
| No. | PhaseA | PhaseB | PhaseC | PhaseA | PhaseB | PhaseC |
| 1 | 3.880% | 2.722% | 3.571% | 0.032% | 0.296% | 0.012% |
| 2 | 3.930% | 6.456% | 2.896% | 0.028% | 0.011% | 0.004% |
| 3 | 4.538% | 9.983% | 53.028% | 0.827% | 1.414% | 6.769% |
| **4** | 32.924% | 23.639% | 42.024% | 0.428% | 1.052% | 3.969% |
| 5 | 0.079% | 1.871% | 1.159% | 0.054% | 0.229% | 0.098% |
| **6** | 1.616% | 2.768% | 5.791% | 0.423% | 1.269% | 0.847% |
| 7 | 1.127% | 4.650% | 2.386% | 0.042% | 0.098% | 0.445% |
| **8** | 8.009% | 2.278% | 4.241% | 0.050% | 0.082% | 0.344% |
| 9 | 0.477% | 1.659% | 2.329% | 0.607% | 2.114% | 0.880% |

### D. Cost Minimization

To study how the SDP OPF dispatches distributed generation, ten DGs were installed on the IEEE-123 node test feeder. Eight of these are three-phase DGs and the other two are single-phase DGs. Table VII gives the DGs' locations and parameters, which were carefully chosen to avoid causing voltage problems.

Table VII. Parameters of Distributed Generators in the 123-Node Feeder

| Location | Three-phase Power Rating | | | DG costs (Quadratic) | |
|---|---|---|---|---|---|
| | A | B | C | $a_1$ | $a_2$ |
| Node | kW | kW | kW | ($/kWh) | ($/kWh$^2$) |
| 23 | 40 | 40 | 40 | 0.02 | 0.0004 |
| 47 | 50 | 50 | 50 | 0.08 | 0.0002 |
| 52 | 30 | 30 | 30 | 0.14 | 0.0005 |
| 62 | 60 | 60 | 60 | 0.04 | 0.0001 |
| 77 | 50 | 50 | 50 | 0.08 | 0.0002 |
| 89 | 30 | 30 | 30 | 0.13 | 0.0003 |
| 101 | 40 | 40 | 40 | 0.05 | 0.0005 |
| 135 | 30 | 30 | 30 | 0.07 | 0.0003 |
| 3 | 0 | 0 | 30 | 0.03 | 0.0005 |
| 9 | 30 | 0 | 0 | 0.12 | 0.0002 |

The price of purchasing power from the substation is set at 0.1$/kWh, while the DGs are assumed to have quadratic cost functions. As shown by Equation (24), the objective is to minimize the total cost of purchasing electricity from substation and the DGs.

$$\begin{aligned} \min:\ & \pi_{grid}(P_{grid}^A + P_{grid}^B + P_{grid}^C) \\ & + \sum_{g\in DG} a_1(P_g^A + P_g^B + P_g^C) + a_2(P_g^A + P_g^B + P_g^C)^2 \end{aligned} \tag{24}$$

$P_g^A$ and $Q_g^A$ represent the real and reactive power outputs of DG $g$ on phase A. The power factor of all DGs is set at 0.9 and the real power output of each DG on each phase is constrained between 0 and its power rating:

$$0 \le P_g^A \le P_g^{A,\max},\ 0 \le P_g^B \le P_g^{B,\max},\ 0 \le P_g^C \le P_g^{C,\max} \tag{25}$$

Three-phase DGs are assumed to trip if their power production is too unbalanced. To reflect this, we impose linear [23] constraints (26-27) on the three-phase DG power outputs:

$$\begin{cases} (1-\beta)\dfrac{P_g^A + P_g^B + P_g^C}{3} \le P_g^A \le (1+\beta)\dfrac{P_g^A + P_g^B + P_g^C}{3} \\ (1-\beta)\dfrac{P_g^A + P_g^B + P_g^C}{3} \le P_g^B \le (1+\beta)\dfrac{P_g^A + P_g^B + P_g^C}{3} \\ (1-\beta)\dfrac{P_g^A + P_g^B + P_g^C}{3} \le P_g^C \le (1+\beta)\dfrac{P_g^A + P_g^B + P_g^C}{3} \end{cases} \tag{26}$$

$$\begin{cases} (1-\beta)\dfrac{Q_g^A + Q_g^B + Q_g^C}{3} \le Q_g^A \le (1+\beta)\dfrac{Q_g^A + Q_g^B + Q_g^C}{3} \\ (1-\beta)\dfrac{Q_g^A + Q_g^B + Q_g^C}{3} \le Q_g^B \le (1+\beta)\dfrac{Q_g^A + Q_g^B + Q_g^C}{3} \\ (1-\beta)\dfrac{Q_g^A + Q_g^B + Q_g^C}{3} \le Q_g^C \le (1+\beta)\dfrac{Q_g^A + Q_g^B + Q_g^C}{3} \end{cases} \tag{27}$$

As suggested in [23], the parameter $\beta$ is set at 0.1856. While the BFM-SDP failed to solve, the symmetrical SDP model produced the optimal solution shown in Table VIII:

Table VIII. Cost Minimization Solution of the Symmetrical SDP

| Phase | A | | B | | C | |
|---|---|---|---|---|---|---|
| | P:kW | Q:kVar | P:kW | Q:kVar | P:kW | Q:kVar |
| DG23 | 39.96 | 19.35 | 28.22 | 13.67 | 35.32 | 17.10 |
| DG47 | 23.87 | 11.56 | 16.73 | 8.10 | 19.85 | 9.61 |
| DG52 | 0.01 | 0.00 | 0.01 | 0.00 | 0.01 | 0.00 |
| DG62 | 60.00 | 29.06 | 60.00 | 29.06 | 60.00 | 29.06 |
| DG77 | 25.12 | 12.17 | 17.37 | 8.41 | 21.15 | 10.24 |
| DG89 | 0.01 | 0.00 | 0.01 | 0.00 | 0.01 | 0.00 |
| DG101 | 21.89 | 10.60 | 15.13 | 7.33 | 18.42 | 8.92 |
| DG135 | 22.04 | 10.67 | 15.36 | 7.44 | 18.41 | 8.92 |
| DG3 | 0.00 | 0.00 | 0.00 | 0.00 | 30.00 | 14.53 |
| DG9 | 0.01 | 0.00 | 0.00 | 0.00 | 0.00 | 0.00 |
| Source | 1267.00 | 461.18 | 811.0 | 258.50 | 986.2 | 285.00 |

From the above table, we observe that DG62 and DG3 produce at their maximum capacity because their marginal costs at full capacity are still lower than that of the grid. DG23, DG47, DG77, DG101 and DG135 are partially dispatched at levels with marginal costs slightly higher than the grid power because their outputs not only supply loads, but also reduce network losses. For these partially dispatched DGs their phase C production is greater than phase B production but less than phase A production. This is because every kW of DG production in phase C reduces losses more than in phase B but less than in phase A.

## V. Voltage Regulation SDP Model

So far, all the objective functions that we have considered are strictly increasing in the power injections. This condition holds when the objective is to minimize losses or minimize the total cost. However, as Zamzam et. al [9] point out, if the objective function takes other forms, the rank-1 constraint (7) may no longer hold and the solution of SDP might become meaningless. Voltage regulation is a good example where SDP relaxations fail. The objective is indeed to maintain an acceptable voltage profile across the network, which does not necessarily strictly mean an increase in power injections. Several authors have used SDP to solve voltage-related OPF problems. Zhang et al. [24] developed an SDP model to tackle the voltage regulation of distribution networks with deep penetration of distributed energy resources. Dall'Anese et al. [25] proposed a model to optimally dispatch inverters of residential photovoltaic systems. Robbins et al. [26] solved for the optimal tap setting of voltage regulators.

Different from the references, our objective is to maintain an acceptable voltage profile, rather than dispatching resources to achieve an optimal voltage profile. An acceptable voltage profile is defined by constraint (4): $\underline{v}_i \le diag(v_i) \le \overline{v}^i$. Our experience suggests that when (4) is binding, the SDP models tend to generate meaningless solutions with rank higher than 1. To overcome this problem, we seek an alternative to replace this constraint.

First, we derive a linear approximation of DG outputs-nodal voltages relationship. As in LPF model [8], we omit the loss term (last term) of Equation (3). The voltage at any node *i* can then be expressed as:

$$v_i = v_0^{\Phi_j} - \sum_{(k,l)\in Path(i)} [S_{kl} z_{kl}^H + z_{kl} S_{kl}^H]^{\Phi_j} \quad (28)$$

where $Path(i)$ denotes the path from the source bus going downstream to bus *i*. The complex power output vector of DG is: $\Lambda_g = [P_g^A + jQ_g^A, P_g^B + jQ_g^B, P_g^C + jQ_g^C]^T$. The outputs of DGs reduces the power transferred from upstream. Considering the DG injections, Equation (28) becomes:

$$v_i - \Delta v_i = v_0^{\Phi_i} - \sum_{(k,l)\in Path(i)} [(S_{kl} - \Delta S_{kl}) z_{kl}^H + z_{kl}(S_{kl} - \Delta S_{kl})^H]^{\Phi_i} \quad (29)$$

Hence: $$\Delta v_i = \sum_{(k,l)\in Path(i)} [\Delta S_{kl} z_{kl}^H + z_{kl} \Delta S_{kl}^H]^{\Phi_i} \quad (30)$$

where the change in $\Delta S_{kl}$ is due to the power injections of the DGs located downstream of line *kl*. We approximate $\Delta S_{kl}$ by summing the outputs of these downstream DGs:

$$\Delta S_{kl} = \sum_{g\in Down(k)} S_g^{\Phi_{kl}} \quad (31)$$

where $Down(k)$ denotes set of DGs that are downstream of node *k*. Assuming the three-phase voltages are balanced, we transform the power output vector $\Lambda_g$ of DG *g* to $S_g$, a 3-by-3 matrix.

$$S_g = \Gamma^{\Phi_g} diag(\Lambda_g) \quad (32)$$

where $$\Gamma = \begin{bmatrix} 1 & a & a^2 \\ a^2 & 1 & a \\ a & a^2 & 1 \end{bmatrix} \text{ and } a = 1\angle 120^\circ \quad (33)$$

Please refer to [8] for a discussion of this approximation. Equations (28) – (33) express DG outputs and their impacts on voltage profiles in a linear manner.

Let turn our attention to the voltage constraint (4). If we ignore all DGs, the OPF problem becomes a power flow problem, which can be solved easily. Let us denote the voltages as calculated by the power flow by $v_i^{(0)}$. Dispatching the DGs should alter the voltage profile in such a way that every nodal voltage satisfies constraint (4), as in (35):

$$\underline{v}_i \le diag(v_i^{(0)} + \Delta v_i) \le \overline{v}^i \quad (34)$$

$$\underline{v}_i \le diag(A v_i^{012,(0)} A^H + \Delta v_i) \le \overline{v}^i \quad (35)$$

We replace Constraint (4) by (34), where $\Delta v_i$ is a linear function of the change in DG power outputs $\Delta\Lambda_g$. With (34) replacing (4) and (35) replacing (14), the voltage regulation problem can be solved with the symmetrical SDP technique. Figure 2 shows how the iterative approach updates the loads and enforces the voltage constraints.

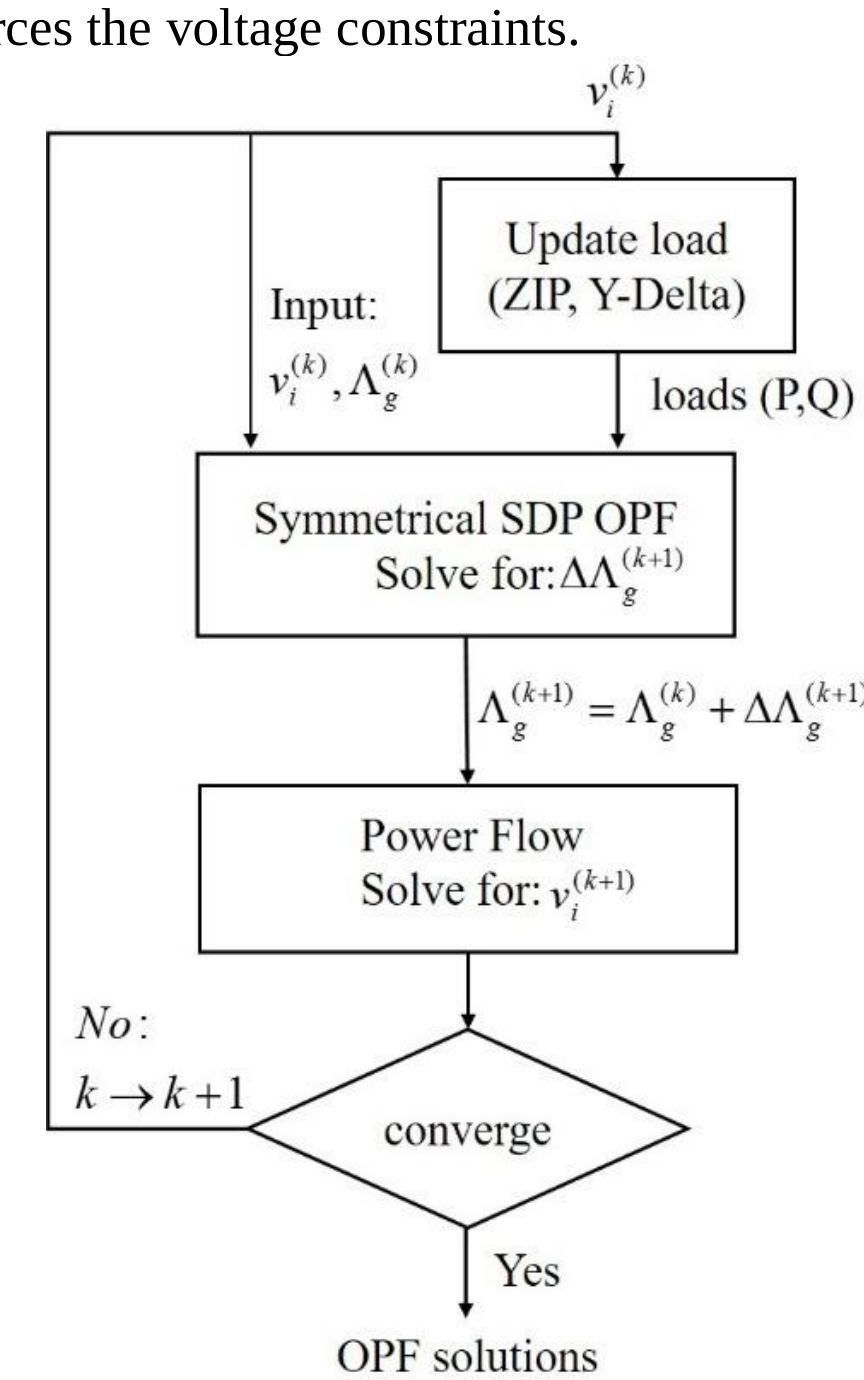

Figure 2. Flow chart of iterative voltage regulation OPF model.

The algorithm contains a series of successive linear approximations:

1. At iteration *k*, the voltages and the DG power outputs $v_i^{(k)}$ are gathered. In constraints (34-35), we replace the voltage parameters with latest results $v_i^{(k)}$ Then we solve the voltage regulation SDP model and record the solutions: the changes in DG outputs $\Delta\Lambda_g^{(k+1)}$.
2. We update the DG outputs $\Lambda_g^{(k+1)} = \Lambda_g^{(k)} + \Delta\Lambda_g^{(k+1)}$ and send the updated outputs to a power flow module. Based on the latest DG outputs $\Lambda_g^{(k+1)}$, the power flow module solves the voltage $v_i^{(k+1)}$.
3. Iterations are repeated until the convergence criterion are met.

## VI. Voltage Regulation SDP model Case Studies

To validate the effectiveness of our proposed algorithm, we perform case studies on the 13- and 34-node test feeders. For the 13-node feeder case, the source bus voltage is set at 1.0 p.u. and the taps of the voltage regulators are set at 10(A), 10(B) and 11(C) to create over-voltage problems. Based on the OpenDSS power flow results, the phase B voltage at several nodes exceeds its upper limit. The voltage at node 675 is the highest at 1.0686p.u.

Three single-phase battery energy storages are installed on the three phases of node 675. Each battery has a power converter rated at 200kVA, with a power factor fixed at 0.95. To make this OPF non-trivial, the price of charging the batteries is set at $\pi_g = \$0.13 / kWh$, higher than the price of electrical energy from the substation. The objective function is to minimize the total cost of supplying power:

$$\min: \pi_{grid}(P_{grid}^A + P_{grid}^B + P_{grid}^C) + \sum_{g \in Storage} \pi_g (P_g^A + P_g^B + P_g^C) \qquad (36)$$

According to the voltage regulation symmetrical SDP OPF, the battery on phase B should be charged at 169.20+j55.61 kW+jkVAr, while the batteries on the other two phases should remain idle because the voltages on these two phases are within limits. Table IX shows that the voltages resulting from the SDP are quite close to those calculated by OpenDSS when the control variables are set at the values calculated by the SDP OPF.

Table IX. Comparison between the Open

| | OpenDSS | | | Voltage Regulation SDP | | |
|---|---|---|---|---|---|---|
| | VA(pu) | VB(pu) | VC(pu) | VA(pu) | VB(pu) | VC(pu) |
| SRC | 0.9999 | 0.9999 | 0.9999 | 1.0000 | 1.0000 | 1.0000 |
| 650 | 0.9998 | 0.9999 | 0.9998 | 0.9999 | 0.9999 | 0.9999 |
| RG60 | 1.0622 | 1.0622 | 1.0684 | 1.0624 | 1.0624 | 1.0686 |
| 632 | 1.0187 | 1.0464 | 1.0219 | 1.0189 | 1.0466 | 1.0221 |
| 645 | \ | 1.0373 | 1.0199 | \ | 1.0374 | 1.0201 |
| 633 | 1.0156 | 1.0445 | 1.0193 | 1.0159 | 1.0447 | 1.0195 |
| 671 | 0.9850 | 1.0498 | 0.9870 | 0.9853 | 1.0500 | 0.9872 |
| 646 | \ | 1.0355 | 1.0178 | \ | 1.0357 | 1.0180 |
| 634 | 0.9916 | 1.0263 | 1.0005 | 0.9918 | 1.0264 | 1.0007 |
| 684 | 0.9831 | \ | 0.9850 | 0.9833 | \ | 0.9852 |
| 680 | 0.9850 | 1.0498 | 0.9870 | 0.9853 | 1.0500 | 0.9872 |
| 692 | 0.9850 | 1.0498 | 0.9870 | 0.9852 | 1.0500 | 0.9872 |
| 652 | 0.9775 | \ | \ | 0.9778 | \ | \ |
| 611 | \ | \ | 0.9830 | \ | \ | 0.9832 |
| 675 | 0.9791 | 1.0498 | 0.9855 | 0.9793 | 1.0499 | 0.9857 |

In addition, the power flows at the feeder head calculated by the symmetrical SDP are very close to the power flow solutions from OpenDSS, with absolute errors on both P and Q around 0.1%.

The second study focuses on the 34-node test feeder. The source bus voltage is set at 1.0 p.u. and the taps at regulator 852 are adjusted to 9(A), 8(B) and 9(C) create serious under-voltage issues. The voltages at nodes 814 and 852 are as low as 0.90p.u. Voltage constraints are not imposed on node 890 because doing this creates over-voltage issues at other nodes. This feeder provides a good case to assess the linear DG-voltage relations (31) – (34) because the losses on this feeder represent about 16% of the total load.

Three DGs are installed at nodes 820, 824 and 860. DG820 is a single-phase DG with capacity of 100kW. DG824 and DG860 are three-phase DGs with per phase capacity of 200kW and 100kW, respectively. The power factor of the DGs are fixed at 0.95. The prices for three DGs are set at 0.12, 0.15 and 0.13 $\$/kWh$. The SDP results are validated using OpenDSS. Table X shows the dispatch of the DGs produced by the SDP OPF and Table XI the nodal voltages from SDP and OpenDSS.

Table X. DG Outputs from Voltage Regulation OPF on 34-Node Feeder Case.

| Phase | A | | B | | C | |
|---|---|---|---|---|---|---|
| | P:kW | Q:kVar | P:kW | Q:kVar | P:kW | Q:kVar |
| DG820 | 100 | 32.87 | \ | \ | \ | \ |
| DG824 | 98.443 | 32.36 | 93.82 | 30.84 | 71.64 | 23.55 |
| DG860 | 100 | 32.87 | 100 | 32.87 | 79.47 | 26.12 |
| Source | 394.74 | 32.26 | 422.09 | 10.28 | 416.9 | -35.51 |

Table XI. Nodal Voltages on 34-Node Test Feeder.

| | OpenDSS | | | Voltage Regulation SDP | | |
|---|---|---|---|---|---|---|
| | VA(pu) | VB(pu) | VC(pu) | VA(pu) | VB(pu) | VC(pu) |
| SRC | 1.0000 | 1.0000 | 1.0000 | 1.0000 | 1.0000 | 1.0000 |
| 800 | 1.0000 | 1.0000 | 1.0000 | 1.0000 | 1.0000 | 1.0000 |
| 802 | 0.9988 | 0.9989 | 0.9990 | 0.9988 | 0.9989 | 0.9990 |
| 806 | 0.9980 | 0.9982 | 0.9984 | 0.9980 | 0.9982 | 0.9984 |
| 808 | 0.9825 | 0.9863 | 0.9872 | 0.9825 | 0.9862 | 0.9873 |
| 810 | \ | 0.9861 | \ | \ | 0.9860 | \ |
| 812 | 0.9643 | 0.9734 | 0.9737 | 0.9645 | 0.9732 | 0.9741 |
| 814 | 0.9498 | 0.9632 | 0.9630 | 0.9500 | 0.9629 | 0.9635 |
| 814R | 1.0211 | 0.9933 | 0.9931 | 1.0212 | 0.9930 | 0.9936 |
| 850 | 1.0210 | 0.9933 | 0.9931 | 1.0212 | 0.9930 | 0.9936 |
| 816 | 1.0208 | 0.9931 | 0.9929 | 1.0210 | 0.9929 | 0.9935 |
| 818 | 1.0204 | \ | \ | 1.0205 | \ | \ |
| 824 | 1.0160 | 0.9872 | 0.9875 | 1.0162 | 0.9869 | 0.9881 |
| 820 | 1.0114 | \ | \ | 1.0116 | \ | \ |
| 826 | \ | 0.9870 | \ | \ | 0.9867 | \ |
| 828 | 1.0155 | 0.9867 | 0.9869 | 1.0156 | 0.9864 | 0.9875 |
| 822 | 1.0084 | \ | \ | 1.0086 | \ | \ |
| 830 | 1.0021 | 0.9736 | 0.9732 | 1.0023 | 0.9733 | 0.9739 |
| 854 | 1.0018 | 0.9733 | 0.9728 | 1.0019 | 0.9730 | 0.9735 |
| 856 | \ | 0.9732 | \ | \ | 0.9729 | \ |
| 852 | 0.9792 | 0.9504 | 0.9491 | 0.9792 | 0.9501 | 0.9500 |
| 852R | 1.0342 | 0.9979 | 1.0025 | 1.0343 | 0.9975 | 1.0034 |
| 832 | 1.0342 | 0.9979 | 1.0025 | 1.0343 | 0.9975 | 1.0034 |
| 858 | 1.0328 | 0.9964 | 1.0010 | 1.0329 | 0.9960 | 1.0019 |
| 888 | 0.9976 | 0.9615 | 0.9669 | 0.9976 | 0.9611 | 0.9679 |
| 864 | 1.0328 | \ | \ | 1.0329 | \ | \ |
| 834 | 1.0313 | 0.9946 | 0.9993 | 1.0313 | 0.9942 | 1.0002 |
| 890 | 0.9144 | 0.8866 | 0.8851 | 0.9143 | 0.8862 | 0.8864 |
| 860 | 1.0312 | 0.9946 | 0.9992 | 1.0313 | 0.9941 | 1.0002 |
| 842 | 1.0312 | 0.9946 | 0.9993 | 1.0313 | 0.9942 | 1.0002 |
| 836 | 1.0310 | 0.9942 | 0.9991 | 1.0311 | 0.9938 | 1.0000 |
| 844 | 1.0310 | 0.9943 | 0.9990 | 1.0311 | 0.9938 | 1.0000 |
| 840 | 1.0310 | 0.9941 | 0.9990 | 1.0311 | 0.9937 | 1.0000 |
| 862 | 1.0310 | 0.9942 | 0.9991 | 1.0311 | 0.9937 | 1.0000 |
| 846 | 1.0313 | 0.9942 | 0.9993 | 1.0314 | 0.9938 | 1.0002 |
| 838 | \ | 0.9940 | \ | \ | 0.9936 | \ |
| 848 | 1.0314 | 0.9942 | 0.9993 | 1.0314 | 0.9938 | 1.0003 |

These results show that DG820 is dispatched at maximum capacity because it is the cheapest way to boost the voltages. DG824 and DG860 are partially dispatched to minimize the cost needed to achieve an acceptable voltage profile. The voltages at nodes 814 and 852 are binding at 0.95p.u. The voltages calculated by the SDP are quite close to the OpenDSS benchmark, with errors no greater than 0.001p.u. According to OpenDSS, the three-phase power flows at the feeder head are 394.2+j34.7 on phase A, 421.9+j9.1 on phase B and 418.2−j32.4 on phase C, which are close to the values calculated by the SDP. This is in contrast to the results of Tables V and VI, which showed that the original SDP OPF achieved the least accurate solution on the 34-node feeder. Ten iterations are required for both the 13-node and 34- node feeders. These case studies demonstrate the proposed technique for regulating voltages using the SDP model is effective and accurate.

## VII. Conclusions and Discussions

This paper proposes two new formulations of an SDP-based model designed to solve OPF problems in unbalanced three-phase distribution systems. Comprehensive case studies demonstrate that these new SDP formulations have better numerical stability and achieve more accurate solutions than existing approaches.

Future research could focus on the following areas: 1) Optimality: in this paper, we demonstrate the numerical stability and accuracy of our models by comparing solutions with OpenDSS power flow results. However, we cannot guarantee that these formulations lead to optimal solutions. 2) Modeling: it would be useful to integrate models of other distribution system components (delta-connected regulators, substation transformers, distribution transformers, triplex line etc.) into the formulation. The iterative process could also be optimized to reduce the number of iterations and speed up the solution process. 3) SDP limitations: it is very important to identify the limitations of the SDP approach and the boundaries beyond which it fails to generate meaningful solutions. 4) Mixed-integer programming: some common distribution system OPF problems (such as volt-var optimization and network reconfiguration) require integer variables to represent the taps of voltage regulators and the on/off status of shunt capacitors and switches. Once a stable and efficient mixed-integer SDP solver becomes available, we could extend the SDP approach to address mixed-integer OPF problems.